\documentclass[12pt,oneside]{article}
\RequirePackage{epsfig}
\RequirePackage{amsmath}
\RequirePackage{amsfonts,amssymb,amsfonts}
\RequirePackage{amsthm}

\usepackage{color}
\usepackage{graphicx}
\usepackage{graphics}
\usepackage{subfigure}
\usepackage{footnote}

\newtheorem{lemma}{Lemma}
\newtheorem{proposition}{Proposition}
\newtheorem{definition}{Definition}

\newtheorem{remark}{Remark}

\newtheorem{example}{Example}
\newenvironment{pf}%
{\par\noindent\textbf{Proof:}~}%
{\eop\par\smallskip\par\noindent}
\newcommand{\ZZ}{\mathbb{Z}}
\newcommand{\RR}{\mathbb{R}}

\newcommand{\CC}{\mathbb{C}}

\def\bgamma{\mbox{\boldmath $\gamma$}}

\newcommand{\even}{\mbox{$\phantom{d}$even}}
\newcommand{\odd}{\mbox{$\phantom{d}$odd}}
\newcommand\ba{\mathbf{a}}

\newcommand\bq{\mathbf{q}}

\newcommand\bbq{\mathbf{q}}
\newcommand\ie{{\it\thinspace i.e.}}

\newenvironment{code}{
                           \mathcode`\:="603A  
                           \def\colon{\mathchar"303A}
                           \par
                           \upshape
                           \begin{list} 
                              {} {\leftmargin = 0.0cm}
                           \item[]
                           \begin{tabbing}
                              \hspace*{.3in} \= \hspace*{.3in} \=
                              \hspace*{.3in} \= \hspace*{.3in} \=
                          \hspace*{.3in} \= \hspace*{.3in} \= \kill
                          }{\end{tabbing}\end{list}}
\newcommand{\be}{\begin{equation}}
\newcommand{\ee}{\end{equation}}
\newcommand{\barr}{\begin{array}}
\newcommand{\earr}{\end{array}}

\newcommand{\comment}[1]{}

 \newcommand{\B}[1]{\mbox{\boldmath $#1$}}

\newcommand{\diag}{\mbox{\rm {diag}}}

\newcommand{\eop}{\hfill$\Box$}

\begin{document}


\newpage

\begin{center}
\textbf{FROM APPROXIMATING TO INTERPOLATORY NON-STATIONARY SUBDIVISION SCHEMES WITH
THE SAME GENERATION PROPERTIES}

\bigskip

\textsc{Costanza Conti}\\
\footnotesize
\textrm{Universit\`{a} di Firenze, Dipartimento di Energetica ``Sergio Stecco''\\
Via Lombroso 6/17, 50134 Firenze, Italia
}\\
\textit{costanza.conti@unifi.it}

\bigskip
\normalsize
\textsc{Luca Gemignani}\\
\footnotesize
\textrm{Universit\`{a} di Pisa, Dipartimento di Matematica\\
Largo Bruno Pontecorvo 5, 56127 Pisa, Italia
}\\
\textit{gemignan@dm.unipi.it}

\bigskip
\normalsize
\textsc{Lucia Romani}\\
\footnotesize
\textrm{Universit\`{a} di Milano-Bicocca, Dipartimento di Matematica e Applicazioni\\
Via R. Cozzi 53, 20125 Milano, Italia}\\
\textit{lucia.romani@unimib.it}

\normalsize

\bigskip

\begin{abstract}
In this paper we describe a general, computationally feasible strategy to deduce a
family of interpolatory non-stationary subdivision schemes from a symmetric
non-stationary, non-interpolatory one satisfying quite mild
assumptions. To achieve this result we extend our previous work
[C.~Conti,  L.~Gemignani, L.~Romani, Linear Algebra Appl.  431  (2009),  no. 10, 1971--1987]
to full generality by removing additional assumptions on the input symbols.
For the so obtained interpolatory schemes we prove that they are capable of reproducing the same exponential polynomial space as the one generated by the original approximating scheme.
Moreover, we specialize the computational methods for the case of symbols obtained
by shifted non-stationary affine combinations of exponential B-splines, that are at the basis of most
non-stationary subdivision schemes. In this case we find that the associated family of interpolatory symbols
can be determined to satisfy a suitable set of generalized interpolating conditions at the set of the zeros
(with reversed signs) of the input symbol. Finally, we discuss some computational examples by showing that the proposed approach can yield novel smooth non-stationary interpolatory subdivision schemes possessing very interesting reproduction properties.



\end{abstract}

\end{center}

\bigskip

\normalsize

{\renewcommand{\thefootnote}{}

\footnotetext{Date: April 16, 2010}
}

\section{Introduction}
Binary interpolatory subdivision schemes are efficient
iterative procedures for the generation of interpolatory curves:
starting with the set of points to be interpolated, at each
recursion step a new point is inserted in between any two given
points so that the limit curve, whenever exists, not only
interpolates the initial set of points but also all the points
generated through the whole process.
Taking into account that a curve is displayed on the screen by
visualizing a discrete set of its points,
from a computational viewpoint
interpolatory subdivision schemes turn out to be more efficient than classical interpolating methods in several situations.
In fact the limit points obtained within five or six subdivision iterations are in general enough for a good
discrete representation of the limit shape.
This is one of the reasons why interpolatory subdivision schemes are widely used in applications and often preferred
to standard methods.

Two important areas where interpolatory subdivision schemes play a crucial role are
Computer Aided Geometric Design (CAGD) and wavelets construction (see \cite{DL02} and \cite{Micchelli96}, respectively).
In these fields a fundamental issue that recently emerged is
concerned with the study of numerical algorithms for converting known approximating schemes into new interpolatory ones.
Starting from the works \cite{MS01} and \cite{R04}, where the conversion is obtained for a specific approximating scheme by means of a \emph{push-back} or a \emph{tweak} operator,
geometric approaches based on the idea that an interpolatory refinement can be interpreted as an averaging step on the control points
followed by a further adjustment of some of them to fit the interpolation constraints were presented \cite{LM07,LLYL08}.
Very recently a completely different technique  relying upon the interplay between   polynomial and structured matrix computations
has been proposed in  \cite{CGR09}.
In that work  for a  given
\emph{symmetric Hurwitz} approximating symbol an associated family of interpolatory symbols  is  determined in such a way to satisfy an auxiliary
polynomial equation.
As it clearly appears, although the latter strategy turns out to be more general than the previous ones, it is limited to
the context of stationary subdivision schemes.
Being non-stationary subdivision schemes more powerful than stationary ones and
very attractive in several applications such as in CAGD (because
of their ability to reproduce conic sections, spirals or widely used trigonometric curves)
it is of fundamental importance to provide a general and efficient method to convert a given non-stationary,
non-interpolatory scheme into a family of interpolatory ones.
To our knowledge, there exists only a new paper \cite{BCR10} addressing this problem, which presents a strategy
that is restricted  to the case of symmetric subdivision masks of odd width, namely symmetric subdivision symbols of even degree.

The goal of this paper is to elaborate on our recent work \cite{CGR09} to progress along different directions.
In particular, (i) we extend the applicability of the proposed
construction, (ii) we investigate the reproduction properties of the so-obtained
interpolatory schemes and (iii) we design algorithms specifically suited for the case of approximating symbols
generated from exponential B-splines, that are at the basis of most
non-stationary subdivision schemes. More specifically,  in this paper  we prove that the strategy
described in  \cite{CGR09} can still be  pursued  under very relaxed conditions
on the approximating symbols we deal with,
say $\{a^{(k)}(z),\ k\ge0 \}$. If, for a given fixed $k\geq 0$,  $a^{(k)}(z),\ a^{(k)}(-z)$ are relatively prime  polynomials,  then a double family of interpolatory symbols associated with $a^{(k)}(z)$
can be generated by solving  two different Bezout-like polynomial equations.  In the symmetric case where $a^{(k)}(z)$ is a symmetric polynomial, the
double family reduces to one single family since the solutions of these two equations are suitably related.  In the Hurwitz case where
$a^{(k)}(z)$ is a Hurwitz polynomial, the distribution of the roots implies the primality condition.
Whenever  such a condition is satisfied for any $k\geq 0$ then  the  correspondence of $a^{(k)}(z)$
with any member of the associated double family   allows one  to  define  a  family  of interpolatory subdivision schemes
derived  from the given non-stationary approximating one.  The computation of the  interpolatory symbol amounts to solve the corresponding
polynomial equation.  If the approximating symbol is specified by spectral information,  as it is  generally the case of
exponential B-splines,   then it is shown that the equation can be efficiently solved by using the tool of (incomplete) partial fraction decomposition.
This gives a representation of the associated interpolatory symbol in terms of a set of generalized interpolating conditions attained at the  zeros
(with reversed signs) of
the approximating  symbol. For the newly generated  interpolating  schemes
we prove an
important reproduction result: the exponential polynomial space reproduced by the interpolatory
scheme is the same function space generated by the approximating one it is originated from.
On the contrary, a general result concerning convergence and/or smoothness
of a non-stationary interpolatory subdivision scheme induced by a non-stationary approximating one is not yet available.
However, in many specific examples we have considered, the analysis can be performed by using ad-hoc techniques.
In this way, by starting with approximating schemes suitably generated by five term affine combinations of exponential B-splines, we are able to find novel smooth
non-stationary interpolatory subdivision schemes possessing very interesting reproduction properties.\\

The paper is organized as follows. In Section \ref{Back} the needed background on non-stationary subdivision schemes is
given. In Subsection \ref{step_sec} we review  and generalize    the basic
strategy  proposed in \cite{CGR09} for the construction of
an interpolatory subdivision mask from a given approximating one.
Effective computational procedures for implementing this strategy
are discussed in Subsection \ref{subsection:roots}. These procedures are the key ingredients of  our
algorithm,   named  {\em Appint} and stated in Subsection \ref{Appint_algo},   to move from a non-stationary approximating
subdivision scheme to a family of
non-stationary interpolatory ones.  The reproduction properties of these schemes are studied in
Section \ref{properties} whereas in Section \ref{ERS}
the application of the algorithm to several instances of
non-stationary approximating subdivision schemes generating exponential polynomials is considered.
Finally, conclusions and further work are drawn in Section \ref{end}.

\section{Background}\label{Back}
In this section we briefly recall some needed background on stationary and non-stationary subdivision schemes.
For more material on subdivision schemes we refer the reader to the seminal work by Cavaretta,
Dahmen and Micchelli \cite{CDM91}, to the more recent survey  by Dyn and Levin \cite{DL02} and to
the well-known book by Warren and Weimer \cite{WW02}.

Subdivision schemes are simple iterative algorithms to efficiently generate curves and surfaces.
Any subdivision scheme is defined by an infinite sequence of coefficients collected in the so
called \emph{refinement masks} $\{{\ba}^{(k)},\ k\ge 0\}$. We
assume that any mask ${\ba}^{(k)}:=\left( a_{i}^{(k)} \in \RR,\ i\in \ZZ\right)$
is of real numbers and has finite support for all $k\ge 0$ \ie\
$a_i^{(k)}=0$ for $i\not\in [-n(k),n(k)]$ for suitable $n(k)\geq 0$.
The \emph{k-level subdivision operator} associated with the $k$-level mask
${\ba}^{(k)}$ is
\begin{equation}\label{soper}
S_{\ba^{(k)}} \ :\ \ell(\ZZ) \rightarrow
\ell(\ZZ) \ ,\qquad \displaystyle{(S_{\ba^{(k)}}\
{\bq})_{i}:=\sum_{j \in\ZZ} a^{(k)}_{i-2j}\ {q}_{j}},\ \ \ i \in
\ZZ \ ,
\end{equation}
where $\ell(\ZZ)$ denotes the linear space of real sequences
indexed by $\ZZ$ whose elements will be denoted by boldface
letter, $\bq:=\left(q_i\in\RR, i\in \ZZ\right)$.
The \emph{subdivision scheme} consists of the
subsequent application of $S_{\ba^{(0)}},\cdots, S_{\ba^{(k)}}$ from a given
starting sequence, say $\bq$, generating the scalar sequences
\begin{equation}\label{sscheme}
{\bq}^{(0)}:={\bq}\ , \ \ \ {\bf \bq}^{(k+1)}:=S_{\ba^{(k)}}\, {\bq}^{(k)}\
\hbox{\ for \ \ } k\ge 0.
\end{equation}
In case the masks
$\{\ba^{(k)},\ k\ge 0\}$ are kept fixed over the iterations, that is $\ba^{(k)}=\ba$
for all $k\geq 0$, the subdivision scheme is said to be \emph{stationary}, otherwise \emph{non-stationary}.\\
Attaching the data $q_i^{(k)}$ generated at the $k$-th step to the parameter values $t^{(k)}_i$ with
$$t^{(k)}_i<t_{i+1}^{(k)}, \quad \hbox{and} \quad t_{i+1}^{(k)}-t_i^{(k)}=2^{-k},\quad k\ge 0$$
(these are usually set as $t_i^{(k)}:=\frac{i}{2^k}$) we see that the subdivision process generates denser and denser sequences of data so that a notion of convergence can be established by taking into account the piecewise
linear function $Q^{(k)}$ that interpolates the data, namely
\[
  Q^{(k)}(t_i^{(k)}) = q_i^{(k)}, \qquad
  Q^{(k)}|_{[t_i^{(k)},t_{i+1}^{(k)}]} \in \Pi_1, \qquad
  i\in\ZZ,\quad k\geq0,
\]
where $\Pi_1$ is the space of linear polynomials. If
the sequence $\{Q^{(k)},\ k\ge 0\}$ converges, then we denote its
limit by
\[
  f_\bbq := \lim_{k\to\infty} Q^{(k)}
\]
and say that $f_\bbq$ is the \emph{limit function} of
the subdivision scheme based on the rule (\ref{sscheme}) for the data $\bbq$.
Several subdivision properties can be read off from the
\emph{symbols}
\[
 a^{(k)} (z) = \sum_{i\in \ZZ}a^{(k)}_i \, z^i,\quad k\ge 0, \qquad z
\in \CC\setminus \{0\}
\]
associated to the masks
$\{{\ba}^{(k)},\ k\ge 0\}$. Also, the corresponding \emph{sub--symbols}
\[
a^{(k)}_{\even}(z)=
\sum_{i\in \ZZ}a^{(k)}_{2i} \, z^i,
\quad
a^{(k)}_{\odd}(z)=
\sum_{i\in \ZZ}a^{(k)}_{2i+1} \, z^i,   \quad z
\in \CC\setminus \{0\},
\]
related to the symbols by the relation
\[
a^{(k)}_{\even}(z^2) + z \cdot a^{(k)}_{\odd}(z^2)
=a^{(k)}(z),
\]
are useful tools for subdivision analysis. Note that since the masks are always supposed to be
finitely supported, all symbols are Laurent polynomials. Nevertheless, for the analysis of
subdivision properties of our concern we can always assume to work with polynomial symbols,
at least after the application of a suitable shift at each iteration.

 A celebrated class of stationary subdivision schemes is given by
degree-$n$ \emph{polynomial B-spline} subdivision schemes, whose (unique) symbol is
\begin{equation}\label{def:B-splinesn}
   B_n(z)=\frac{(1+z)^{n+1}}{2^n}\,, \quad k\ge 0.
\end{equation}
The non-stationary counterpart of (\ref{def:B-splinesn}) is the symbol of the so-called \emph{exponential B-splines}.
They are piecewise functions whose pieces are exponential polynomials (the latter ones will be recalled
in the next definition). These are defined in terms of a linear differential operator and turn out to be of great interest in
geometric modeling  for the design of important analytical shapes like
conic sections, spirals and classical trigonometric curves.

\begin{definition}(Space of exponential polynomials)\label{def:Vspace}
Let $T\in \ZZ_+$ and $\bgamma=(\gamma_0,\gamma_1,\cdots,\gamma_T)$
with $\gamma_T\neq 0
 $ a finite set of real or imaginary numbers and let $D^n$
the $n$-th order differentiation operator. The space of
exponential polynomials $V_{T,\bgamma}$ is the subspace
\begin{equation}\label{def:VMgamma}
V_{T,\bgamma}:=\{f:\RR\rightarrow \CC, f\in C^T(\RR):\quad
\sum_{j=0}^T \gamma_jD^j\,f =0\}.
\end{equation}
\end{definition}

A characterization of the space $V_{T,\bgamma}$ is provided by the following:

\begin{lemma}\label{Prop:characterization}\cite{BR89}
Let $\gamma(z)=\sum_{j=0}^T \gamma_jz^j$ and denote by
$\{\theta_\ell,\tau_\ell\}_{\ell=1,\cdots,N}$ the set of zeros
with multiplicity of $\gamma(z)$ satisfying
$$\gamma^{(r)}(\theta_\ell)=0,\quad r=0,\cdots,\tau_\ell-1,\quad \ell=1,\cdots,N.$$ It results
$$
T=\sum_{\ell=1}^N \tau_\ell,\qquad
V_{T,\bgamma}:=Span\{x^{r}e^{\theta_\ell\,x},\
r=0,\cdots,\tau_\ell-1,\ \ \ell=1,\cdots,N\}.
$$
\end{lemma}

As proved in \cite{MWW01} (see also \cite{WW02}) exponential B-splines can be generated via a
non-stationary subdivision scheme based on the  symbols
\begin{equation}\label{def:ER}
B_{n}^{(k)}(z)=2\prod_{\ell=1}^N\left(
\frac{e^{\frac{\theta_\ell}{2^{k+1}}}z+1}{e^{\frac{\theta_\ell}{2^{k+1}}}+1}\right)^{\tau_\ell},\quad
k\ge 0\,.
\end{equation}
Its limit function belongs to the subclass of $C^{T-2}$
degree-$n$ L-splines \cite{S81} (with $n=T-1$) whose pieces are exponentials
of the space $V_{T,\bgamma}$.
Notice that, when $\theta_{1}=0$ with $\tau_1=n+1$, then
$B_n^{(k)}(z)$ in (\ref{def:ER}) does not depend on $k$ being the
symbol of a degree-$n$ B-spline given in (\ref{def:B-splinesn}).
An important aspect of subdivision schemes is their convergence capability to specific classes of functions.
In particular, a subdivision scheme is said to possess the property of \emph{generating exponential polynomials} if,
for any initial data uniformly sampled from some exponential polynomial function, the scheme yields a function belonging to
the same space in the limit. Even more, the subdivision scheme is \emph{reproducing exponential polynomials} if,
for any initial data uniformly sampled from some exponential polynomial function, the scheme yields the \emph{same} function in the limit.
To this purpose, we recall the following two important definitions (see, for example, \cite{CR10} and \cite{VonBluUnser07}).

\begin{definition} [$V_{T,\bgamma}$-Generation] \label{def:ERrgenerationlimit}
Let $\{a^{(k)}(z),\ k\ge 0\}$ be a set of subdivision symbols. The subdivision scheme associated with the set of symbols $\{a^{(k)}(z),\ k\ge 0\}$ is
said to be \emph{$V_{T,\bgamma}$-generating} if it is convergent and
for $f\in V_{T,\bgamma}$ and for the initial sequence ${\bf f}^{0}:=\{f(t^{0}_i),\ i\in \ZZ\}$,  it results
$$\lim_{k\rightarrow \infty}S_{\ba^{(k)}}\cdots S_{\ba^{(0)}}{\bf f}^{0}=\tilde f\,,\quad \tilde f\in V_{T,\bgamma}\,.
$$
\end{definition}

\begin{definition} [$V_{T,\bgamma}$-Reproduction] \label{def:ERreproductionlimit}
Let $\{a^{(k)}(z),\ k\ge 0\}$ be a set of subdivision symbols. The subdivision scheme associated with the symbols $\{a^k(z),\ k\ge 0\}$ is
said to be \emph{$V_{T,\bgamma}$-reproducing} if it is convergent and
for $f\in V_{T,\bgamma}$ and for the initial sequence ${\bf f}^{0}:=\{f(t^{0}_i),\ i\in \ZZ\}$,  it results
$$\lim_{k\rightarrow \infty}S_{\ba^{(k)}}\cdots S_{\ba^{(0)}}{\bf f}^{0}=f\,.
$$
\end{definition}

Since the space of exponential polynomials trivially includes standard polynomials,
Definitions \ref{def:ERrgenerationlimit} and \ref{def:ERreproductionlimit} include, as special cases, the notion of polynomial generation
and polynomial reproduction, respectively. For a complete analysis of the latter concepts in the stationary situation
--which are very much related to the approximation order of the subdivision scheme-- the interested reader can  see \cite{DHSS08}.

We conclude  by recalling that a subdivision
scheme is said to be interpolatory if the refinement masks
$\{{\ba}^{(k)},\ k\ge 0\}$ satisfy
\begin{equation}\label{interc1}
a^{(k)}_{2i}=\delta_{i,0},\quad
\hbox{or equivalently,} \quad
a^{(k)}_{\even}(z)=1,\quad
   k\ge0,
\end{equation}
meaning that all points generated by the subdivision process at a
given level $k$ will be kept in the next level $k+1$.
We also mention that from \eqref{interc1} it follows that
a mask $\ba^{(k)}$ is interpolatory if and only if all its symbols
$a^{(k)}(z)$ satisfy the algebraic condition
\begin{equation}\label{def:interp}
    a^{(k)}(z)+a^{(k)}(-z)=2,\quad
   \forall k\ge0.
\end{equation}

\section{\bf From approximating to interpolatory subdivision schemes}\label{stepwise}

In this section we introduce the key ingredients of our proposed algorithm  termed  {\em Appint} to
generate a family of non-stationary interpolatory subdivision schemes starting from an
initial non-stationary approximating one.  At the core of this algorithm there is a procedure which, for a given fixed
non-interpolatory subdivision symbol  $a^{(k)}(z)$, $k\geq 0$, effectively constructs a
corresponding interpolatory symbol denoted by $m^{(k)}(z)$.
The procedure is applied step-by-step for $k=0, 1,\ldots$.
For the sake of notational simplicity  we  can therefore omit the superscript $k$
by denoting $a^{(k)}(z)=a(z)$ and  $m^{(k)}(z)=m(z)$.
The construction stems from a  theoretical result  presented in
 \cite[Theorem 2]{CGR09} which describes the conditions being satisfied for the associated interpolatory
symbol $m(z)$. In Subsection
\ref{step_sec} this   result is
reviewed  and generalized  to some extent by removing unnecessary restrictions on the  input symbol $a(z)$.
In the case where  $a(z)$ is of the form \eqref{def:ER} and it  is known in factorized form  by means of the set of zeros
$\{\theta_\ell,\tau_\ell\}_{\ell=1,\cdots,N}$, then  an efficient
method for computing a suitable representation
of $m(z)$ is described in Subsection
\ref{subsection:roots}. Finally,  by  putting  all these ingredients together,
{\em Appint} is formally stated in
Subsection \ref{Appint_algo}.

\subsection{From approximating to interpolatory subdivision symbols}\label{step_sec}

In the matrix environment   the  linear operator $S_\ba$  defined in \eqref{soper} and associated with
the symbol  $a(z)=\sum_{i\in \ZZ}a_i z^i$, $z\in \CC\setminus \{0\}$
is represented by a bi-infinite Toeplitz-like matrix
$S_\ba=(a_{i-2j})$, $i,j \in \mathbb Z$. Since $a(z)$ is a Laurent polynomial,
say $a(z)=\sum_{j=-\kappa}^{\kappa} a_jz^j$,
$\max\{|a_{-\kappa}|,|a_{\kappa}|\}>0$, it follows that $S_\ba$ is banded with bandwidth
$\lceil \frac{\kappa}{2}\rceil$ at most. Let
$p(z)=\sum_{j=-h}^h p_jz^j$, $ \max\{|p_{-h}|,|p_{h}|\} >0$, be another Laurent polynomial and denote by
${\mathcal P}$ the bi-infinite Toeplitz matrix associated with $p(z)$, namely, ${\mathcal P}=(p_{i-j})$.
Observe that ${\mathcal P}$ is again banded with bandwidth $h$. For the product operator
\[
\mathcal S \colon = \mathcal P \cdot S_\ba =(s_{i,j}), \quad i,j \in \mathbb Z,
\]
we have
\[
s_{i,j}=\sum_{r=i-h}^{i+h} p_{i-r} \, a_{r-2j}=\sum_{\ell=-h}^{h} p_{\ell} \, a_{i-2j-\ell}=s_{i+2, j+1},
\quad i, j\in \mathbb Z.
\]
This means that the product operator $\mathcal S$  is a  bi-infinite Toeplitz-like  matrix of the same form as
the subdivision operator $S_\ba$  with entries $s_{i,j}=s_{i-2j}$, $i, j\in \mathbb Z$.
By setting
\[
q(z)=a(z) \cdot p(z)= \sum_{j=-h-\kappa}^{h+\kappa} q_j z^{j}, \quad (q_j=0 \ \mbox{if} \ |j|>h+\kappa),
\]
we find that
\[
q_j=\sum_{i=-h}^h p_i \, a_{j-i}, \quad -(h+\kappa)\leq j\leq h+\kappa,
\]
and, therefore,
\[
q_{i-2j}=s_{i,j}=s_{i-2j}, \quad i, j \in \mathbb Z.
\]
There follows that the  product operator $\mathcal S$ can be seen as the subdivision operator associated with the
Laurent polynomial $q(z)$, i.e.,
\[
\mathcal S= S_\bq, \quad q(z)=a(z)\cdot p(z),
\]
where  $a(z)$ is the symbol of $S_\ba$  and $p(z)$ can be suitably chosen in such a way to
satisfy the interpolation condition.
By expressing $q(z)$  in terms of its sub--symbols
\[
 q(z)=   q_{\even}(z^2) + z \cdot q_{\odd}(z^2)     \qquad z
\in \CC\setminus \{0\},
\]
we find that
\[
q(z)+ q(-z)=2 \cdot  q_{\even}(z^2).
\]
Then by imposing the interpolation condition   \eqref{interc1}, i.e., $q_{\even}(z)=1$,  we arrive at the  relation
\begin{equation}\label{bezout}
a(z) \cdot p(z) + a(-z) \cdot p(-z) =2
\end{equation}
which is a generalized Bezout equation providing necessary and sufficient conditions  for a
Laurent polynomial $p(z)$ to convert  the subdivision operator associated with $a(z)$
into the interpolating subdivision operator generated by $q(z)=a(z) \cdot p(z)$.

Suitable coefficient-wise representations of $p(z)$ are introduced
to  investigate conditions under which the (generalized) Bezout equation is solvable as well as to develop
 effective computational methods for
its solution. Observe that if $p(z)$ is of the form
\begin{equation}\label{polt}
p(z)=p_{\kappa}z^{\kappa} + p_{\kappa+1}z^{\kappa+1}+\ldots + p_{\kappa+m}z^{\kappa+m},
\end{equation}
with $m=2 \kappa -1$,
and, moreover,  it satisfies
\begin{equation}\label{bez_gen_f}
a(z) \cdot p(z) + (-1)^j a(- z) \cdot p(- z)  =2 z^{j}, \quad  0\leq j\leq 2m+1,
\end{equation}
then  $z^{-j} p(z)$  solves   \eqref{bezout}.  Computing polynomial solutions of
\eqref{bez_gen_f} of the form \eqref{polt}  reduces in
 a matrix setting to solving a  structured linear system whose coefficient matrix is
Sylvester-like.  Let $\B a_0=\left[a_{-\kappa}, \ldots, a_0, \ldots, a_{\kappa}\right]^T$ $\in \mathbb R^{2\kappa +1}$
denote the coefficient vector
of the Laurent polynomial $a(z)$.  The associated extended
coefficient vector $\B {\widehat a}_+\in \mathbb R^{2\kappa +m+1}$ is  defined by
$\B {\widehat a}_+^T=\left[\B a_0^T, 0, \ldots, 0\right]$.  Similarly let us introduce the extended coefficient vector
$\B {\widehat a}_-\in \mathbb R^{2\kappa +m+1}$  associated with the polynomial $a(-z)$.
Moreover let $Z=(z_{i,j})\in \mathbb R^{2 (m+1)\times 2 (m+1)}$
be the down-shift matrix  given by
$z_{i,j}=\delta_{i-1,j}$, where $\delta_{i,j}$ is the Kronecker delta symbol.
Set $\mathcal R_+ \in \mathbb R^{2(m+1) \times (m+1)}$
the striped Toeplitz matrix
\[
\mathcal R_+=\left[\B {\widehat a}_+| Z\B {\widehat a}_+|\ldots |Z^m\B {\widehat a}_+\right],
\]
and,  similarly,  define
\[
\mathcal R_-=\left[\B {\widehat a}_-| Z\B {\widehat a}_-|\ldots |Z^m\B {\widehat a}_-\right].
\]
The coefficient matrix of the linear system  \eqref{bez_gen_f} is
$\mathcal R^+=\left[\mathcal R_+| \mathcal R_-\right]\in \mathbb R^{2(m+1)\times 2(m+1)}$
or $\mathcal R^-=\left[\mathcal R_+| -\mathcal R_-\right]\in \mathbb R^{2(m+1)\times 2(m+1)}$
depending on the parity of $j$.
It is well known that $\mathcal R^+$  and $\mathcal R^-$ are resultant matrices and, therefore,
they are   invertible if and only
if $a(z)$ and $a(-z)$ are relatively prime polynomials.

Due to the special structure of the polynomial pair
$(a(z),a(-z))$   it is shown that both linear systems can be reduced to smaller systems of half the size.
Let  $P_{m+1} \in \mathbb  R^{2 (m+1)\times 2 (m+1)}$,  $P_{m+1}=(\delta_{i, \sigma(j)})$  be the permutation matrix associated with the
``perfect shuffle'' permutation given by
$$\sigma \ : \
\{1,\ldots, 2m +2\} \rightarrow \{1,\ldots, 2m+2\},\quad
\sigma(j)=\left\{%
\begin{array}{ll}
    (j+1)/2 + m+1 , &  \hbox{\ \ if $j$ is odd;} \\
    \\
    j/2, &  \hbox{\ \ if $j$ is
even.} \\
\end{array}%
\right. $$
Furthermore,  let   $G_{m+1} \in \mathbb R^{2 k\times 2 k}$ be the matrix
defined by
\[
G_{m+1} =\left(\begin{array}{c|c}  I_{m+1}   & -D_{m+1}\\\hline D_{m+1} & I_{m+1}
\end{array} \right),
\]
where $D_{m+1}=\diag[-1, (-1)^2, \ldots, (-1)^{k-1}, (-1)^{m+1}].$
There follows that
\begin{equation}\label{h1}
P_{m+1} \cdot \mathcal R^- \cdot G_{m+1}^{-1}=\mathcal H^- \oplus \mathcal H,
\end{equation}
where $\mathcal H \in \mathbb  R^{(m+1) \times (m+1)}$ is a certain matrix  and
\[
\mathcal H^-=\left[\begin{array}{cccccc}  a_{-\kappa+1} & a_{-\kappa} & 0 & \ldots & \ldots    & \ldots
\\ a_{-\kappa+3} & a_{-\kappa+2}  & a_{-\kappa+1} & a_{-\kappa} & 0 & \ldots \\
 a_{-\kappa+5} &a_{-\kappa+4} & a_{-\kappa+3}  & \ldots  &  & \ldots \\
\vdots & \vdots & \vdots & \vdots &   &  \ldots \\
\vdots & \vdots & \vdots & \vdots &   & \ldots\\
a_{-\kappa+2m+1} & a_{-\kappa+2m}  & a_{-\kappa+2m-1} &\ldots & & \ldots
\end{array}
\right].
\]
Similarly we find that
\begin{equation}\label{h2}
P_{m+1} \cdot \mathcal R^+ \cdot G_{m+1}^{-1}=\widehat{\mathcal H}  \ {\widehat \oplus} \  \mathcal H^+,
\end{equation}
where $\widehat{\mathcal H} \in \mathbb  R^{(m+1) \times (m+1)}$ is a certain matrix,
$\widehat \oplus$ denotes the direct sum with respect to the main anti-diagonal,
and, moreover,
\[
\mathcal H^+=\left[\begin{array}{cccccc}  a_{-\kappa} & 0 & \ldots &  \ldots  &    & \ldots
\\ a_{-\kappa +2} & a_{-\kappa+ 1}  & a_{-\kappa} & 0 &  & \ldots \\
 a_{-\kappa +4} &a_{-\kappa+3} & a_{-\kappa+2}  & \ldots  &  & \ldots \\
\vdots & \vdots & \vdots & \vdots &   &  \ldots \\
\vdots & \vdots & \vdots & \vdots &   & \ldots\\
a_{-\kappa+2m} & a_{-\kappa+2m -1}
   & a_{-\kappa+2m-2} &\ldots & & \ldots
\end{array}
\right].
\]
In this way we  arrive at the following generalization  of \cite[Theorem 2]{CGR09}.
\begin{proposition}\label{cappero}
Let $\widehat a(z)=z^\kappa a(z)$ be a degree-$n$
polynomial, $n=m+1$,  relatively prime with $\widehat a(-z)$.
Then $\mathcal H^-$  and $\mathcal H^+$ are  invertible and, moreover, the
polynomial $p_i^\star(z)$, $\star \in \{+,-\}$,  with coefficients given by the entries of the
$i$-th column  of $({\mathcal H^\star})^{-1}$, $1\leq i \leq n$, is the unique
polynomial of degree less than $n$ such that
\begin{equation}\label{relFOND}
\widehat a(z)p_i^\star(z)\  \star \  \widehat a(-z)p_i^\star (-z)= 2\, z^{2i-\ell^\star}, \qquad 1\leq i \leq n, \qquad \star \in \{+,-\}
\end{equation}
where
\[
\ell^\star= \left\{\begin{array}{ll} 2 \ {\rm if}  \ \star=+; \\
1, {\rm elsewhere}. \end{array}\right.
\]
\end{proposition}

As an immediate consequence of Proposition \ref{cappero}  we obtain   the following.

\begin{proposition}\label{Prop:Interpolation}
Given a degree-$n$ polynomial $\widehat a(z)$ relatively prime
with $\widehat a(-z)$ and such that $\widehat a(1)=2$, $\widehat a(-1)=0$,  then the Laurent
polynomials
\begin{equation}\label{def:m^i_n}
m_i^\star(z):=\frac{\widehat a(z)p_i^\star (z)}{z^{2i-1}}, \quad 1\le i\le n,
\end{equation}
 where $p_i^\star (z)$ solves \eqref{relFOND}, $\star \in \{+,-\}$,  are  the
associated interpolatory symbols
and satisfy
$$m_i^\star (1)=2,\quad m_i^\star(-1)=0, \quad 1\le i\le n.$$
\end{proposition}

\begin{remark}\label{in1}
It is worth noting that Proposition \ref{cappero} defines a double family of associated interpolatory symbols
depending on the sign of $\star$. In the symmetric case where $\widehat a(z)$ is a symmetric polynomial,
that is, $a_j=a_{-j}$, $0\leq j\leq \kappa$, the number of  associated symbols  halves  since
all the matrices $\mathcal H$, $\widehat{\mathcal H}$,
$\mathcal H^+$ and $\mathcal H^-$ are suitably related and, in particular, $\mathcal H^+$ can be obtained from
$\mathcal H^-$ by reversion of  rows and columns.
\end{remark}

These  results  provide a practical way to construct a family of finitely supported interpolatory masks
from a given approximating one consisting in computing the matrix $({\mathcal H^\star})^{-1}$
 and reading its entries.  This approach seems to be  especially  tailored
for symmetric Hurwitz subdivision symbols
which result into  computations with totally positive (TP)
Hurwitz matrices. The procedures described in \cite{Pena}
can be adjusted for the efficient and stable computations
of the coefficients of the interpolatory masks generated in
the B-spline case and ``shifted'' affine combinations of them (see \cite[Section 4]{CGR09}).
However, in the case of exponential B-splines and their  affine combinations the
approximating symbol is generally known by assigning  the spectrum of the symbol, that is, its  zeros with
their multiplicity. It is therefore interesting to design a completely different machinery for
solving \eqref{relFOND} using the information on the roots.

\subsection{A root-based polynomial equation solver}\label{subsection:roots}
Let us suppose that
\[
\widehat a(z)=\widehat a_0 + \widehat a_1z + \ldots +\widehat a_{n}z^{n}= \widehat a_{n}\prod_{j=0}^m
(z-z_j)^{k_j},
\]
with $z_i\neq z_j$ if $i\neq j$ and $k_0 + \ldots + k_m=n$. Then  it is shown that
the   unique
solution  $p_i(z)$  of \eqref{relFOND} can be obtained by imposing certain
interpolation conditions at the zeros of $\widehat a(z)$ and  $\widehat a(-z)$.

Let us start by recalling the concept of {\em
Hermite-Lagrange interpolation polynomial}   of a given
differentiable function $f(z)$ on the set of nodes $\eta_0, \ldots, \eta_\ell$ with multiplicities $h_0,
\ldots, h_\ell$, $h_0+ \ldots +h_\ell=r+1$, respectively. Suppose that the function
$f(z)$ possesses derivatives $f^{(j)}(\eta_i)$, $0\leq j\leq h_i-1$, $0\leq i\leq \ell$.
Then there exists a unique  polynomial
$H_f(z)$  of degree at most $r$ satisfying the interpolation
conditions
\[
H_f^{(j)}(\eta_i)=f^{(j)}(\eta_i), \quad 0\leq j\leq h_i-1, \
0\leq i\leq \ell.
\]
This polynomial is generally  referred to as the  Hermite-Lagrange
interpolation polynomial of $f(z)$ on the prescribed set of nodes.
By setting $\omega(z):=(z-\eta_0)^{h_0} \cdots
(z-\eta_\ell)^{h_\ell}$ we find the Lagrange-type representation
\[
H_f(z)=\sum_{i=0}^\ell \sum_{j=0}^{h_i-1} \sum_{h=0}^{h_i-j-1}
f^{(j)}(\eta_i) \frac{1}{h!j!}
\left(\begin{array}{c}\displaystyle\frac{(z-\eta_i)^{h_i}}
{\omega(z)}\end{array}\right)_{z=\eta_i}^{(h)}
\displaystyle\frac{\omega(z)}{(z-\eta_i)^{h_i-j-h}}
\]
and, equivalently, the partial-fraction representation
\[
H_f(z)=\omega(z)\sum_{i=0}^\ell
\sum_{s=1}^{h_i}\displaystyle\frac{1}{(z-\eta_i)^s}\left(\sum_{j=0}^{h_i-s}
\mathcal S(h_i-j-s, j, i)\right)=
\omega(z)\sum_{i=0}^\ell
\sum_{s=1}^{h_i}\displaystyle\frac{c_{i,h_i-s}}{(z-\eta_i)^s},
\]
where
\[
\mathcal S(h, j, i)=f^{(j)}(\eta_i)\frac{1}{h!j!}
\left(\begin{array}{c}\displaystyle\frac{1}
{\omega_i(z)}\end{array}\right)_{z=\eta_i}^{(h)}, \quad
w_i(z)=\displaystyle\frac{\omega(z)}{(z-\eta_i)^{k_i}},
\]
and, moreover,  by  Leibniz's rule
\[
c_{i,j}=\sum_{\ell=0}^{j}
\mathcal S(j-\ell, \ell, i)=\frac{1}{j!}\left(\begin{array}{c}\displaystyle\frac{H_f(z)}
{\omega_i(z)}\end{array}\right)_{z=\eta_i}^{(j)}.
\]
Let $\ell=2m+1$  and $\eta_0=z_0, \ldots, \eta_{(\ell-1)/2}=z_m$,
$\eta_{(\ell+1)/2}=-z_0, \ldots, \eta_\ell=-z_m$ with
multiplicities $h_0=h_{(\ell+1)/2}=k_0, \ldots,
h_\ell=h_{(\ell-1)/2}=k_m$. Observe that
\[
r + 1=h_0+ \ldots + h_\ell=2 k_0 + \ldots + 2k_m=2n
\]
and
\[
\omega(z)=\prod_{i=0}^m(z-z_i)^{k_i} \prod_{i=0}^m (z+z_i)^{k_i}=
\widehat a_{n}^{-2}(-1)^{n} \widehat a(z) \widehat a(-z).
\]
By replacing  the right hand side  $f(z)=2 z^{2t-\ell^\star}$  of  \eqref{relFOND},
 where $t$ is fixed and $1\leq t\leq n$,  with its
Hermite-Lagrange form  we find that
\[
(-1)^{n} \widehat a_{n}^2 \left(\frac{p_t^\star(z)}{\widehat a(-z)}  \ \star  \
\frac{p_t^\star(-z)}{\widehat a(z)}\right)= \sum_{i=0}^\ell
\sum_{s=1}^{h_i}\displaystyle\frac{c_{i,h_i-s}}{(z-\eta_i)^s}, \qquad \star \in \{+,-\}.
\]
Since $\widehat a(z)$ and $\widehat a(-z)$ are relatively prime we can  separate
the partial fraction decompositions of the two rational functions on the
left-hand side. This gives the following

\begin{proposition}\label{Prop:rappresentazione1}
Let $\widehat a(z)=\widehat a_{n}\prod_{j=0 }^m(z-z_j)^{k_j}$
be a polynomial of degree $n$,  where  $z_i\neq z_j$ if $i\neq j$, $k_0 + \ldots + k_m=n$ and
$\widehat a(z)$ and $\widehat a(-z)$  are relatively prime.
Then, the unique polynomial solution $p_t^\star(z)$, $1\leq t\leq n$, $\star \in \{+,-\}$, of \eqref{relFOND}  satisfies
\[
p_t(z)= (-1)^{\ell^\star}\widehat a_{n}^{-1} \prod_{j=0}^m(z+z_j)^{k_j} \sum_{i=0}^{m}
\sum_{s=1}^{k_i}\displaystyle\frac{(-1)^s c_{i, k_i-s}}{(z+z_i)^s},
\]
where
\[
c_{i,j}=\frac{1}{j!}\left(\begin{array}{c}\displaystyle\frac{2 z^{2t-\ell^\star}}
{\omega_i(z)}\end{array}\right)_{z=z_i}^{(j)}, \quad 0\leq j\leq k_i-1, \ 0\leq i\leq m,
\]
 $\omega(z)$ is the monic polynomial associated with $\widehat a(z)\widehat a(-z)$ and $w_i(z)=\displaystyle\frac{\omega(z)}{(z-z_i)^{k_i}}$.
\end{proposition}

\begin{example}\label{in2}
To illustrate  the computational meaning of the previous result let us consider
the interpolatory symbols associated with the cubic exponential B-spline with $k$-level symbol
$$B_3^{(k)}(z)=\frac{1}{2} (z+1)^2 \frac{z^2+2v^{(k)}z+1}{2(v^{(k)}+1)},$$
where  the  parameter $v^{(k)} \in (0,
+\infty)$ is   defined through the expression
$$v^{(k)}=\frac{1}{2} \left(e^{\theta/2^{k+1}}+e^{-\theta/2^{k+1}} \right)$$
with $\theta \in \{\theta_{\ell}, \ \ell=1,...,N\}$, as in Lemma \ref{Prop:characterization}.
As shown in \cite{BCR07a} this means that $B_3^{(k)}(z)$
  corresponds  to (\ref{def:ER}) with $N=3$, $\theta_{1}=0$, $\theta_2=t$, $\theta_3=-t$ and $\tau_1=2$, $\tau_2=\tau_3=1$, and,
moreover,  once assigned the starting value $v^{(-1)}\in (-1, +\infty)$,
the parameter $v^{(k)}$ can be recursively updated at each successive iteration through the formula
\begin{equation}\label{recurrence}
v^{(k)}=\sqrt{\frac{v^{(k-1)}+1}{2}},\ k\ge 0.
\end{equation}
For any fixed $k\geq 0$, the symmetric interpolatory scheme of smallest support
associated with $B_3^{(k)}(z)$ is obtained  from the choice $i=2$ and $\star=-$ in \eqref{relFOND}.  By using
Proposition  \ref{Prop:rappresentazione1} we find that the corresponding solution $p_2^k(z)$ is given by
\[
p_2^{(k)}(z)=\frac{(1-z)^2}{2 v^k({v^k}-1)} -\frac{z^2-2 v^k z +1}{2({v^k}-1)}=\frac{1}{2v^k}\left(-z^2 + 2(v^k+1)z -1\right),
\]
which from Proposition \ref{Prop:Interpolation} defines the interpolatory symbol
\[
m_{3,2}^{(k)}(z)\colon =\frac{B_3^{(k)}(z)p_2^{(k)}(z)}{z^3}, \quad k\geq 0.
\]
\end{example}

The partial fraction decomposition is not a flexible computational tool and several difficulties arise in order to find efficient
updating procedures for computing the solutions of  \eqref{relFOND}  associated with slightly  modified symbols (as usually it is the case in non-stationary subdivision schemes depending on a parameter, see Section \ref{ERS}).
In this respect the tool of
{\em incomplete partial fraction decomposition}  \cite{Henrici} is much more suited.   The general strategy proceeds as follows.
From the partial fraction decomposition we get two polynomials $h(z)$ and $k(z)$ of degree less than $n$ such that
$\displaystyle \frac{1}{\widehat a(z) \widehat a(-z)}=\displaystyle\frac{h(z)}{\widehat a(z)} + \displaystyle\frac{k(z)}{\widehat a(-z)}.$
Since $\widehat a(z)$ is  given in factored form  we can  determine
 $k(z)$ as the  Hermite-Lagrange polynomial interpolating the function $g(z)=1/\widehat a(z)$ on the zeros of $\widehat a(-z)$.  Then the polynomial
$p_t^\star(z)$ which solves   \eqref{relFOND} can be obtained by means of the polynomial division  between
$2 z^{2t-\ell^\star}  k(z)$ and $\widehat a(-z)$. Again this operation reduces to computing the
Hermite-Lagrange polynomial interpolating $2 z^{2t-\ell^\star} \cdot k(z)$ on the zeros of $\widehat a(-z)$. In the case where
the initial symbol $\widehat a(z)$ is modified by a linear or a quadratic  factor,  both the two steps in the above procedure can be
modified accordingly. For instance the polynomial $k(z)$ can be specified in the form $k(z)=k_1(z) + \widehat a(-z) \psi(z)$, where
$k_1(z)$ is  the Hermite-Lagrange polynomial interpolating the function $g(z)=1/\widehat a(z)$ on the zeros of $\widehat a(-z)$ and
$\psi(z)$ is a linear factor whose coefficients are  determined so that $k(z)$ satisfies the modified equation.
This approach has been  implemented and used  for computing the interpolatory symbols associated with certain affine combinations of
exponential B-splines. Some computational results  are shown in  Section \ref{ERS}.

\subsection{The Appint algorithm for the non-stationary case}\label{Appint_algo}

So far we have introduced a quite general strategy for deriving  a family of interpolatory symbols from a given
approximating symbol based on the solution
of equation (\ref{relFOND}).
In the non-stationary setting, we compute a family
of non-stationary interpolatory subdivision schemes associated with a
non-stationary approximating one via the solution
of (\ref{relFOND}) at each recursion step.
Therefore, the procedure we consider turns out to be as follows:
assuming  $\{\widehat a^{(k)}(z),\ k\ge 0\}$ are the
degree-$n(k)$ symbols of an approximating non-stationary scheme with $\widehat a^{(k)}(z)$ and $\widehat a^{(k)}(-z)$
relatively prime for all $k\ge 0$, we construct the non-stationary
interpolatory subdivision scheme based on the symbols
$\{m^{(k)} _{i(k)}(z),\ k\ge 0\}$  where, for each $k$, $m^{(k)}_{i(k)}(z),\ 1\le i(k) \le
n(k)$, is one of the interpolatory symbols satisfying
(\ref{relFOND}).  Here and hereafter  for the sake of simplicity we omit the superscript $\star \in \{+,-\}$ since we
assume that the sequence  $(i(k), \star)$, $k \geq 0$,  is given in input and, therefore,
$m^{(k)} _{i(k)}(z)$ denotes  the unique solution  of (\ref{relFOND}) for the given  pair
 $(i(k), \star)$.
Surely, the performance of the non-stationary
subdivision scheme will depend on the selection of the sequence $(i(k), \star)$, $k \geq 0$.
The computational kernel consists of finding the solution of (\ref{relFOND}) for the
input symbol $\widehat a^{(k)}(z)$ and the fixed pair $(i(k), \star)$.  This task can be accomplished by
the inversion of the corresponding matrices $\mathcal H^\star$ or, alternatively, by  means of the  procedure described
in the previous section based on computing the incomplete partial fraction decomposition. The auxiliary routine {\em Solve}
takes in input a suitable representation of $\widehat a^{(k)}(z)$ together with the pair  $(i(k), \star)$ and returns as output the
corresponding solution  $p^{(k)}_{i(k)}(z)$ of (\ref{relFOND}).
For clarity we describe the overall procedure in algorithmic form.

\framebox{\parbox{7.0cm}{
\begin{code}
\bigskip
\centerline{\bf Appint Algorithm}\\

\bigskip
\medskip \noindent  \sf \ \ \ \ \ \ \ \ \ \ Input:
$\{\widehat a^{(k)}(z),\ k\ge 0\}$, degree-$n(k)$ symbols;
\\
\noindent  \sf \ \ \ \ \ \ \ \ \ \ \ \ \ \ \  \quad \
$\{(i(k),  \star),\ k\ge 0\}$,  with
$1\le i(k) \le n(k)$
\\

\bigskip \noindent  \sf \ \ \ \ \ \ \ \ \ \ For $k=0,1,\dots$
\\

\medskip \noindent  \sf \ \ \ \ \ \ \ \ \ \ \ \ \ \ \ Check whether  $\widehat a^{(k)}(z)$ is relatively prime with $\widehat a^{(k)}(-z)$
\\

\medskip \noindent \sf \ \ \ \ \ \ \ \ \ \ \ \ \ \ \
Set  $p^{(k)}_{i(k)}(z)\colon = {\bf Solve}[\widehat a^{(k)}(z), \ (i(k),  \star)]$
\\

\medskip \noindent  \sf \ \ \ \ \ \ \ \ \ \ \ \ \ \ \
Construct the interpolatory symbol  $m^{(k)}_{i(k)}(z):=\frac{\widehat a^{(k)}(z) p_{i(k)}^{(k)}(z)}{z^{2 \, i(k)-1}}$ \hspace{0.8cm}
\\

\bigskip \noindent  \sf \ \ \ \ \ \ \ \ \ \
Output: $\{m^{(k)}_{i(k)}(z),\ k \geq 0\}$

\end{code}
}}

\medskip

Some theoretical properties of the  computed sequence
$\{m^{(k)} _{i(k)}(z),\ k\ge 0\}$ are discussed in  Section \ref{properties} whereas
computational examples are reported in  Section \ref{ERS}.

\section{Properties of non-stationary interpolatory subdivision
schemes derived from their approximating counterparts}
\label{properties}

For the family of non-stationary interpolatory subdivision
schemes generated by symbols $\{m^{(k)}_i(z), k\geq 0\}$,
$1\le i\le n(k)$, we can prove an
important reproduction result: the exponential polynomial space reproduced by the interpolatory
scheme is the same function space generated by the approximating scheme it is originated from.
To prove it, we first need a preliminary result given in \cite{DLL03}.
Within the rest of this section $V_{T,\bgamma}$ is the space given in Definition \ref{def:Vspace} and
$z_\ell^{(k)}:=e^{-\frac{\theta_\ell}{2^{k+1}}},\ \ \ell=1,\cdots, N,\quad k\ge 0$.

\begin{proposition}\label{theo:zeros}
Let $\{m^{(k)}(z),\ k\ge 0\}$ be a sequence of interpolatory symbols. The
subdivision scheme associated with such a sequence reproduces
$V_{T,\bgamma}$ if and only if for each $k\ge 0$
\begin{equation}\label{condzeros}
\begin{array}{ll}
m^{(k)}(z_\ell^{(k)})=2, \qquad  m^{(k)}(-z_\ell^{(k)})=0, \quad \ell=1,\cdots, N\\
 \\
 \frac{d^{r}}{dz^{r}} \, m^{(k)}(\pm z_\ell^{(k)})=0,\quad r=1,\cdots, \tau_\ell-1,\ \ \ell=1,\cdots, N.
 \end{array}
\end{equation}
\end{proposition}

\bigskip \noindent We are now in a position to state the reproduction result.

\begin{proposition}\label{prop_pg}
Let $\{\widehat a^{(k)}(z),\ k\ge 0\}$ be a sequence of  symbols with $\widehat a^{(k)}(z)$ relatively prime with $\widehat a^{(k)}(-z)$ for all $k\ge 0$.
If the non-stationary approximating subdivision scheme based on the symbols
$\{\widehat a^{(k)}(z),\ k\ge 0\}$ generates the space $V_{T,\bgamma}$, then for all $1\le i\le n(k)$
the non-stationary interpolatory subdivision scheme based on the symbols
\[
m^{(k)}_i(z)=\frac{\widehat a^{k}(z)p_i^{k}(z)}{z^{2i-1}}, \quad k\geq 0,
\]
whenever convergent, reproduces the same space $V_{T,\bgamma}$.
\end{proposition}

\pf
Due to \cite[Theorem 1]{VonBluUnser07} the symbols $\widehat a^{(k)}(z)$ satisfy
\[
\widehat a^{(k)}(-z_\ell^{(k)})=0,\quad
\frac{d^{r}}{dz^{r}} \, \widehat a^{(k)}(-z_\ell^{(k)})=0,\quad
r=1,\cdots, \tau_\ell-1,\ \ell=1,\cdots, N.
\]
By the Leibnitz's differentiation rule, we easily get an
analogous relation to be satisfied by all $m^{(k)}_i(z)$ (for any $1\le i\le n(k)$) that is
\[
m_i^{(k)}(-z_\ell^{(k)})=0,\quad
\frac{d^{r}}{dz^{r}} \, m_i^{(k)}(-z_\ell^{(k)})=0,\quad
r=1,\cdots, \tau_\ell-1,\ \ell=1,\cdots, N.
\]
It remains to consider the behavior of $m_i^{(k)}(z)$ and its
derivatives at the points $z_\ell^{(k)}$. Now,
since  for each $k$
\[
m^{(k)}_i(z) + m^{(k)}_i(-z)=2, \quad 1\leq i \leq n(k),
\]
it follows that
\[
m_i^{(k)}(z_\ell^{(k)})=2
\]
as well as
\[
\frac{d^{r}}{dz^{r}}m_i^{(k)}(z_\ell^{(k)})=
(-1)^{r+1}\frac{d^{r}}{dz^{r}}m_i^{(k)}(-z_\ell^{(k)})=0, \
r=1,\cdots, \tau_\ell-1,\ \ell=1,\cdots, N.
\]
The use of Proposition \ref{theo:zeros} concludes the proof.
\eop

\begin{remark}
We notice that, if an interpolatory subdivision scheme is $V_{T,\bgamma}$-generating, then
due to the interpolatory nature (that is due to the fulfillment of equation (\ref{def:interp})), it is also $V_{T,\bgamma}$-reproducing.
\end{remark}

\begin{remark}
Unfortunately, contrary to the result in Proposition \ref{prop_pg}, a general result concerning convergence and/or smoothness
of a non-stationary interpolatory subdivision scheme induced by a non-stationary approximating one is not available.
However, in all specific examples discussed in Section \ref{ERS} and many others we tested, convergence and smoothness
analysis of the induced non-stationary interpolatory subdivision schemes is provided.
From the examples we see that the smoothness order of the interpolatory scheme is the half of that of the approximating
one it is originated from. This observation gives us a hint for a theoretical result to be investigated in future researches.
\end{remark}

\section{Interpolatory exponential reproducing non-sta- tionary subdivision schemes}\label{ERS}
Aim of this section is to show the application of our strategy to a family of approximating schemes depending on free parameters.
This leads to  a parameter-dependent family of corresponding interpolatory schemes that can be used to design interesting new
non-stationary interpolatory schemes.
In particular, we show that by means of a five term affine combination of exponential B-splines,
we can generate novel smooth
non-stationary interpolatory subdivision schemes possessing very interesting reproduction properties.

Let us consider the interpolatory scheme based on the symbols $m_{3,2}^{(k)}(z)$ introduced in Example \ref{in2}.
The $C^2$ approximating scheme with symbols $\{B_3^{(k)}(z),\ k\ge 0\}$
was originally introduced in \cite{MWW01} where the authors also showed its capability generation of
the function space $V_{4,\bgamma}=\{1,x,e^{tx},e^{-tx}\}$ (see also \cite{WW02}).
According to the results in Section \ref{properties}, the associated interpolatory scheme  turns out to be the $C^1$
4-point interpolatory scheme reproducing the function space $V_{4,\bgamma}=\{1,x,e^{tx},e^{-tx}\}$
(see also \cite{R09}).  The  reproduction properties of this scheme can be improved
by considering the family of approximating subdivision schemes given by a 5-term affine
combination of $B_3^{(k)}(z)$ of the form
\[
\begin{array}{ll}
\widehat a^{(k)}(z)&= B_3^{(k)}(z)
\left(\alpha^{(k)}+\beta^{(k)} z+(1-2 \alpha^{(k)}-2 \beta^{(k)}) z^2+\beta^{(k)} z^3+\alpha^{(k)} z^4 \right)
\\
 & =
 B_3^{(k)}(z)
\resizebox{11.6cm}{!}{$
\left(\alpha^{(k)}  + \frac{\beta^{(k)} + \sqrt{(4 \alpha^{(k)} + \beta^{(k)})^2 - 4 \alpha^{(k)}}}{2} z +
\alpha^{(k)} z^2 \right) \left( 1 +  \frac{2(1 - 2\beta^{(k)}-4\alpha^{(k)})}{\beta^{(k)} + \sqrt{(4 \alpha^{(k)} + \beta^{(k)})^2 - 4 \alpha^{(k)}}}z+ z^2 \right)
$}
\end{array}
\]
where $\alpha^{(k)}, \beta^{(k)} \in \RR$ are free parameters.
By imposing the primality conditions  for $\widehat a^{(k)}(z), \widehat a^{(k)}(-z)$ it turns out that \eqref{relFOND}
 can be solved  whenever $\alpha^{(k)}\neq 0$ and
$\beta^{(k)} \not \in\{0, \frac{1}{2}-2\alpha^{(k)}, \frac{4(v^{(k)})^2 \alpha^{(k)}-4\alpha^{(k)}+1}{2(1-v^{(k)})}\}$.
 In the case $\alpha^{(k)}=0$  the equation can be degree-reduced in such a way that a polynomial solution can still be found.

By applying the procedure described in Subsection \ref{subsection:roots} we have computed the polynomial $p^{(k)}(z)$   corresponding with the pair
$(i(k), \star)=(4, -)$, $k\geq 0$,  and  set
$$
m^{(k)}(z)=\widehat a^{(k)}(z) p^{(k)}(z) z^{-7}.
$$
By accurately choosing the free parameters $\alpha^{(k)}$ and $\beta^{(k)}$, we can obtain an
interpolatory scheme $m^{(k)}(z)$ that improves the properties of the interpolatory scheme $m_{3,2}^{(k)}(z)$
associated with the combined symbol $B_3^{(k)}(z)$.
Improvements can concern with its reproduction capabilities and/or its smoothness order. In particular:

\begin{enumerate}

\item When $\alpha^{(k)}=0$ and $\beta^{(k)}=\frac{1}{4}$, $\widehat a^{(k)}(z)=\frac{(z+1)^4(z^2+2v^{(k)}z+1)}{16(v^{(k)}+1)}$,
namely it is the $C^4$ exponential B-spline that generates $V_{6,\bgamma}=\{1,x,x^2,x^3,e^{tx},e^{-tx}\}$.
The symbol $m^{(k)}(z)$ is the $C^2$ interpolatory 6-point scheme that reproduces the same space (as previously shown in \cite{R09}).

\item When $\alpha^{(k)}=0$ and $\beta^{(k)}=\frac{1}{4(v^{(k)})^2}$, then
\[
\widehat a^{(k)}(z)=\frac{(z+1)^2(z^2+2v^{(k)}z+1)(z^2+2(2(v^{(k)})^2-1)z+1)}{16(v^{(k)})^2(v^{(k)}+1)},
\]
namely it is the $C^4$ exponential B-spline that generates $V_{6,\bgamma}=\{1,x,e^{tx},e^{-tx},$ $e^{2tx},e^{-2tx}\}$,
while $m^{(k)}(z)$ is the $C^2$ interpolatory 6-point scheme that reproduces the same space (see, again, \cite{R09}).

\item When $\alpha^{(k)}=0$ and $\beta^{(k)}=\frac{1}{2(1+v^{(k)})}$, then
\[
\widehat a^{(k)}(z)=\frac{(z+1)^2(z^2+2v^{(k)}z+1)^2}{8(v^{(k)}+1)^2},
\]
 namely it is the
$C^4$ exponential B-spline that generates $V_{6,\bgamma}=\{1,x,e^{tx},e^{-tx},$ $xe^{tx},xe^{-tx}\}$
and $m^{(k)}(z)$ is the $C^2$ interpolatory 6-point scheme that reproduces the same space \cite{R09}.

\item When $\alpha^{(k)}=\frac{1}{8 (v^{(k)})^2 (v^{(k)}+1) (2v^{(k)}-1)^2}$ and $\beta^{(k)}=\frac{4(v^{(k)})^2-2v^{(k)}-1}{4(v^{(k)})^2(2v^{(k)}-1)^2}$,
then
\[\widehat a^{(k)}(z)=
\resizebox{12cm}{!}{$
\frac{(z+1)^2 (z^2+2v^{(k)}z+1) (z^2+2(4(v^{(k)})^3-3v^{(k)})z+1)(z^2+2(2(v^{(k)})^2-1)z+1)}{32(v^{(k)})^2 (v^{(k)}+1)^2 (2v^{(k)}-1)^2}
$}
,
\]
 \ie\ it is the $C^6$ exponential B-spline generating $V_{8,\bgamma}=\{1,x,e^{tx},e^{-tx},e^{2tx},e^{-2tx},$ $e^{3tx},e^{-3tx}\}$, while $m^{(k)}(z)$ defines the $C^3$ interpolatory 8-point scheme that reproduces the same space (see \cite{CR09}).

\medskip
\item When $\alpha^{(k)}=\frac{1}{8(v^{(k)})^2(v^{(k)}+1)}$ and $\beta^{(k)}=\frac{2v^{(k)}-1}{4(v^{(k)})^2}$,
we deal with  the $C^6$ exponential B-spline
\[
\widehat a^{(k)}(z)=
\frac{(z+1)^2 (z^2+2v^{(k)}z+1)^2 (z^2+2(2(v^{(k)})^2-1)z+1)}{32(v^{(k)})^2(v^{(k)}+1)^2}
\]
generating the function space $V_{8,\bgamma}=\{1,x,e^{tx},e^{-tx},e^{2tx},e^{-2tx},xe^{tx},xe^{-tx}\}$. The symbols $m^{(k)}(z)$
define a  $C^3$ interpolatory 8-point scheme that reproduces the same space (see Proposition \ref{prop_pg}). The smoothness
of the subdivision scheme  $\{m^{(k)}(z),\ k\ge 0\}$ can be obtained through asymptotical equivalence \cite{DL95} with the $C^3$ Dubuc-Deslauriers 8-point interpolatory scheme \cite{DD89,D86}.

\end{enumerate}

The last non-stationary interpolatory subdivision scheme corresponds to a new proposal never presented in the literature.
Other interesting proposals can be obtained by assigning different suitable values to the free parameters $\alpha^{(k)}$ and $\beta^{(k)}$.
In all these kinds of interpolatory schemes, by making the parameter $v^{(-1)}$ local, namely by assuming a different parameter $v_i^{(-1)}$
in correspondence of each edge $\overline{q_i \, q_{i+1}}$ of the starting polyline, we can combine the two important issues of local shape control and
special functions reproduction.
This means that, in the same limit curve, we can include an alternation of exponential polynomial pieces in those regions
where the starting samples belong to one of these curves and smooth limit segments with local tension otherwise.
Also, due to the recurrence relation (\ref{recurrence}), the shape parameter $v^{(k)}$ turns out to be independent of
the parametric values ${\bf t}^{(k)}$, thus reducing computational costs of the algorithm.
In addition to the general reasons discussed in the introduction, these properties contribute to make these
interpolatory subdivision schemes more convenient with respect to the corresponding classical interpolatory methods.

\section{Conclusions and future work}\label{end}

A novel approach has been presented for the computation of
a family of interpolatory non-stationary subdivision schemes from a
non-stationary, non-interpolatory one. The approach reduces
the updating problem  either to the inversion of certain structured
matrices (which can be of Hurwitz type or Sylvester resultant matrices) or to the solution of certain Bezout-like polynomial equations.
If the approximating symbols are defined in terms of spectral information it is shown that the
partial fraction decomposition provides an effective  tool for solving these equations
by yielding a representation of the associated interpolatory symbols
in terms of generalized interpolating conditions.
The newly constructed interpolatory schemes are capable of reproducing the same exponential polynomial space as the one generated by the original approximating scheme.
Although a general result concerning the relationship between convergence and/or smoothness orders of the approximating and interpolatory schemes is not yet available,
ad hoc techniques can be used by showing that in many cases the proposed approach leads to  novel smooth
non-stationary interpolatory subdivision schemes possessing very interesting reproduction properties.
The analysis of more general convergence properties of the subdivision schemes generated by our techniques
is an ongoing research.

\medskip


\end{document}